\documentclass[11pt]{article}
\usepackage[a4paper]{anysize}\marginsize{3.5cm}{3.5cm}{1.3cm}{2cm}
\pdfpagewidth=\paperwidth \pdfpageheight=\paperheight
\usepackage{amsfonts,amssymb,amsthm,amsmath,eucal}
\usepackage{pgf}
\usepackage{bbm,array}
\usepackage{tikz} 
\usepackage{float}
\restylefloat{table}
\usepackage{subfigure}
\usepackage{caption}
\usetikzlibrary{arrows}
\usepackage{multicol}
\usepackage{multirow}

\theoremstyle{plain}
\newtheorem{thm}{Theorem}[section]

\newtheorem{corollary}[thm]{Corollary}

\theoremstyle{definition}

\newtheorem{problem}[thm]{Problem}

\newtheorem{thevarthm}[thm]{\varthmname}

\newenvironment{varthm*}[1]{\trivlist\item[]{\bf #1.}\it}{\endtrivlist}

\renewcommand\geq{\geqslant}

\newcommand\be{\begin{eqnarray*}}
\newcommand\ee{\end{eqnarray*}}

\renewcommand\P{\mathbb P}

\newcommand\newop[2]{\def#1{\mathop{\rm #2}\nolimits}}
\newop\mod{mod}
\newop\log{log}
\newop\ord{ord}
\newop\Gal{Gal}
\newop\SL{SL}
\newop\GL{GL}
\newop\Bl{Bl}
\newop\mult{mult}
\newop\mass{mass}
\newop\div{div}
\newop\codim{codim}
\newop\sing{sing}
\newop\vdim{vdim}
\newop\edim{edim}
\newop\Ass{Ass}
\newop\size{size}
\newop\reg{reg}
\newop\areg{areg}
\newop\asreg{asreg}
\newop\satdeg{satdeg}
\newop\supp{supp}
\newop\gin{gin}
\newop\ini{in}
\newop\vol{vol}
\newop\sat{sat}
\newop\length{length}
\newop\depth{depth}
\newop\characteristic{char}

\def\keywordname{{\bfseries Keywords}}%
\def\keywords#1{\par\addvspace\medskipamount{\rightskip=0pt plus1cm
\def\and{\ifhmode\unskip\nobreak\fi\ $\cdot$
}\noindent\keywordname\enspace\ignorespaces#1\par}}
\def\subclassname{{\bfseries Mathematics Subject Classification
(2010)}\enspace}
\def\subclass#1{\par\addvspace\medskipamount{\rightskip=0pt plus1cm
\def\and{\ifhmode\unskip\nobreak\fi\ $\cdot$
}\noindent\subclassname\ignorespaces#1\par}}

\definecolor{ttqqqq}{rgb}{0.,0.,0.}
\definecolor{zzttqq}{gray}{0.4}

\begin{document}

\author{\L .~Farnik,  J.~Kabat, M.~Lampa-Baczy\'nska, H.~Tutaj-Gasi\'nska}
\title{On the parameter spaces of some B\"or\"oczky configurations}
\date{\today}
\maketitle
\thispagestyle{empty}

\begin{abstract}
    B\"or\"oczky configurations of lines have been recently considered in connection with the problem of the containment between symbolic and ordinary powers of ideals. Here
 we describe parameter families of B\"or\"oczky configurations of $13$, $14$, $16$, $18$ and $24$ lines and investigate rational points of these parameter spaces.

\subclass{52C30,  14N20, 05B30}
\end{abstract}


\section{Introduction}\label{intro}
B\"or\"oczky configurations were originally introduced in connection with the orchard problem. These configurations were described by B\"or\"oczky but the
results were never published, see Crowe and McKee (\cite{CrMcKee1968}) who study these configurations in relation to the celebrated Sylvester-Gallai Theorem.
The same setting appears in \cite{FuPa1984} and \cite{GreTao13} where   configurations with large number of ordinary lines are considered. 

The interest was
recently renewed with reference to the containment problem studied in commutative algebra. The containment relations between symbolic and ordinary powers of
homogeneous ideals have been a quite popular direction of research in algebraic geometry. Let us recall that for a homogeneous ideal  $I$ in
$\mathbb{K}[\mathbb{P}^N]=\mathbb{K}[x_0,\ldots,x_N]$ the $m$-th symbolic power of $I$ is defined as
    $$I^{(m)}=\mathbb{K}[\mathbb{P}^N]\cap\bigcap_{Q\in \text{Ass}(I)}(I^m)_Q,$$
where the intersection is taken over the associated primes of $I$.
The Zariski-Nagata theorem says that
    if $\mathbb{K}$ is algebraically closed and if $I$ is radical then
    $$ I^{(m)}=\Big\{f\in I : \frac{\partial^{|\alpha|}f}{\partial x_i^{\alpha}}=0 \textrm{ on the zeros of } I, |\alpha|<m\Big\}.$$

 In 2001--2002 it was proved (in \cite{ELS} in characteristic
$0$ and in \cite{HoHu} in finite characteristic) that the containment
$$I^{(m)}\subset I^r$$
always holds for $m \geq N \cdot r$, where $N$ is the dimension of the ambient space. On the other hand, it is easy to show an example of the non-containment
$$I^{(2)}\nsubseteq I^2.$$
Thus, a natural question arised: does always the containment
\begin{equation} \label{contain}
I^{(3)}\subset I^2
\end{equation}
hold for an ideal of points on the projective plane? The question is connected with a more general problem, namely, can the power $Nr$ in  the containment  $I^{(Nr)}\subset I^r$ be made smaller?
It turned out that in general it is not the case.
The first counterexample for containment (\ref{contain}) was given by
Dumnicki,  Szemberg and  Tutaj-Gasi\'nska in \cite{DST13}. From then on, in short time there appeared further counterexamples
(see details in \cite{SzSzp2015}), all based on some configurations of lines over the complex numbers. However the problem
was still open over the reals.

The first real counterexample (i.e., a counterexample where the coordinates of all points are real numbers)
 comes from \cite{CGMLLPS2015}. It is also based on a configuration of lines, namely the
B\"or\"oczky configuration of $12$ lines. 
In \cite{HarSec13}  Harbourne and  Seceleanu constructed some new examples in $\mathbb{P}^2$, but in finite
characteristic. Finally, the first rational counterexample, which is a modification of the original B\"or\"oczky configuration of $12$ lines, is  studied by  Dumnicki, Harbourne, Nagel,  Seceleanu, Szemberg, and Tutaj-Gasi\'nska in \cite{res}. 

In recent paper \cite{MJ2015}  Lampa-Baczy\'nska and  Szpond describe in details this rational counterexample, denoted there by $\mathbb{B}12$ (here we will use the notation with  lower index $\mathbb{B}_{n}$ for
the modifications, and the notation $\mathcal{B}_{n}$ for the  original B\"or\"oczky configuration). They also show that no configurations $\mathbb{B}_{15}$
of $15$ lines can be defined over the rational numbers. 

Now the set of triple points of the configurations $\mathbb{B}_{12}$ is not the only known rational counterexample to the containment~(\ref{contain}), another one is considered by the authors in \cite{FKL-BT-G2}.

B\"or\"oczky configurations are proved to be counterexamples for the containment~(\ref{contain})  for $n\in \{12,13 \ldots 24\}$ lines, moreover we believe
that for higher values of $n$ it should also be the case but a computer-aided proof does not finish. It is a natural question to investigate other B\"or\"oczky
configurations to check if there are more rational counterexamples among them. 

In this paper we investigate the families of B\"or\"oczky configurations
$\mathbb{B}_{13}$, $\mathbb{B}_{14}$, $\mathbb{B}_{16}$, $\mathbb{B}_{18}$, $\mathbb{B}_{24}$  consisting of respectively $13$, $14$, $16$, $18$ and $24$
lines, and find  parameter spaces for these families. Thus we complete the study of B\"or\"oczky configurations for previously ``missing'' values of $13$ and
$14$ lines, moreover we discuss the consecutive case of $16$ line, and cases of $18$ and $24$ lines (the number $18$ has the same prime numbers as factors as
$12$, and $24$ is a multiple of $12$). In each case this parameter space is a curve, denoted by $C_{13}, C_{14}, C_{16}, C_{18}, C_{24}$ respectively. For
$C_{13}, C_{14}, C_{16}, C_{18}$ we have checked that these curves have no  rational points which give non-degenerate configurations. For $C_{24}$ the problem
is open (however, the computations with Magma suggest that there are no such rational points).

The paper is organized as follows. The second Section presents the original B\"or\"oczky construction. The third   Section contains a detailed description of
the method of finding $C_{13}$, as well as the description of the degenerations of $\mathbb{B}_{13}$ for some points on $C_{13}$. The fourth Section presents
 the results obtained for $\mathbb{B}_{14}, \mathbb{B}_{16}$, $\mathbb{B}_{18}$ and $\mathbb{B}_{24}$. In each case we describe the way of finding the
parametrization curve $C_{n}$ and what is known about its rational points. We finish the paper with some remarks and open problems.

\section{Original B\"or\"oczky construction}\label{OriginalBoroczky}
The B\"or\"oczky configurations $\mathcal{B}_{n}$ were described in  Example $2$ of \cite{FuPa1984}. Following this example we present here the construction.

Consider a $2n$-gon inscribed in a circle. Let us fix one of the $2n$ points and denote it by $Q_0$. By $Q_{\alpha}$ we denote the point arising by the rotation of $Q_0$ around the center of a circle by the angle $\alpha$.
Then we take the following set of lines
$$\mathcal{B}_{n}=\left\{Q_{\alpha}Q_{\pi - 2\alpha}, \textrm{where } \alpha= \frac{2k \pi}{n} \textrm{ for } k=0, \dots, n-1\right\}.$$
If $\alpha \equiv (\pi - 2\alpha)(\mathrm{mod} \ 2\pi)$ then the line $Q_{\alpha}Q_{\pi - 2\alpha}$ is the tangent to the circle at the point $Q_{\alpha}.$ The
configuration $\mathcal{B}_{n}$ has $\big\lfloor \frac{n(n-3)}{6}\big\rfloor+1$ triple points.

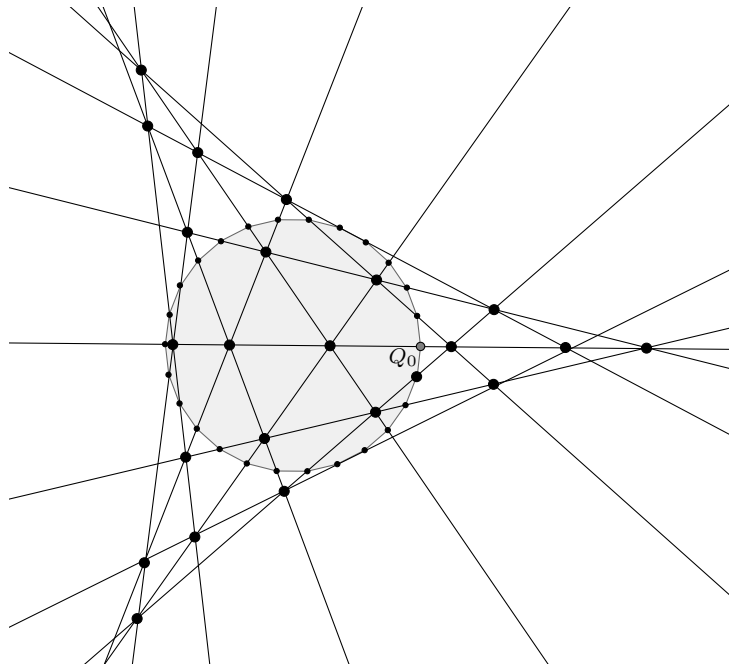
\begin{figure}[h]
\centering
\begin{tikzpicture}[line cap=round,line join=round,x=1.0cm,y=1.0cm, scale=0.15]
\clip(-26.0,-28.0) rectangle (37.0,30.0);
\fill[color=zzttqq,fill=zzttqq,fill opacity=0.1] (10.,0.) -- (9.695,2.68) -- (8.8,5.2) -- (7.2,7.4) -- (5.2,9.2) -- (2.9,10.5) -- (0.2,11.1) -- (-2.5,11.2) -- (-5.1,10.6) -- (-7.5,9.3) -- (-9.5,7.6) -- (-11.056163018827187,5.357977073411441) -- (-12.03155703208026,2.8432143571874904) -- (-12.37678575805597,0.16809900887796264) -- (-12.07178575805597,-2.511900991122036) -- (-11.13428252345037,-5.041033784216154) -- (-9.618760334548844,-7.272315365493406) -- (-7.613295832271324,-9.076071759828515) -- (-5.234439323923323,-10.34747520197045) -- (-2.62044130297844,-11.012636345225395) -- (0.07678216708221663,-11.032898441709754) -- (2.700478262187136,-10.407083932025582) -- (5.098167302036074,-9.171562880732537) -- (7.130504314878069,-7.398139280395476) -- (8.679377260771215,-5.189878064533485) -- (9.654771274024288,-2.6751153483095353) -- cycle;
\draw [color=zzttqq] (10.,0.)-- (9.695,2.68);
\draw [color=zzttqq] (9.695,2.68)-- (8.7574967653944,5.209132793094115);
\draw [color=zzttqq] (8.7574967653944,5.209132793094115)-- (7.241974576492872,7.440414374371368);
\draw [color=zzttqq] (7.241974576492872,7.440414374371368)-- (5.236510074215354,9.244170768706478);
\draw [color=zzttqq] (5.236510074215354,9.244170768706478)-- (2.857653565867352,10.51557421084841);
\draw [color=zzttqq] (2.857653565867352,10.51557421084841)-- (0.24365554492246644,11.180735354103355);
\draw [color=zzttqq] (0.24365554492246644,11.180735354103355)-- (-2.4535679251381928,11.200997450587714);
\draw [color=zzttqq] (-2.4535679251381928,11.200997450587714)-- (-5.0772640202431045,10.575182940903545);
\draw [color=zzttqq] (-5.0772640202431045,10.575182940903545)-- (-7.474953060092044,9.339661889610497);
\draw [color=zzttqq] (-7.474953060092044,9.339661889610497)-- (-9.50729007293404,7.5662382892734374);
\draw [color=zzttqq] (-9.50729007293404,7.5662382892734374)-- (-11.056163018827187,5.357977073411441);
\draw [color=zzttqq] (-11.056163018827187,5.357977073411441)-- (-12.03155703208026,2.8432143571874904);
\draw [color=zzttqq] (-12.03155703208026,2.8432143571874904)-- (-12.37678575805597,0.16809900887796264);
\draw [color=zzttqq] (-12.37678575805597,0.16809900887796264)-- (-12.07178575805597,-2.511900991122036);
\draw [color=zzttqq] (-12.07178575805597,-2.511900991122036)-- (-11.13428252345037,-5.041033784216154);
\draw [color=zzttqq] (-11.13428252345037,-5.041033784216154)-- (-9.618760334548844,-7.272315365493406);
\draw [color=zzttqq] (-9.618760334548844,-7.272315365493406)-- (-7.613295832271324,-9.076071759828515);
\draw [color=zzttqq] (-7.613295832271324,-9.076071759828515)-- (-5.234439323923323,-10.34747520197045);
\draw [color=zzttqq] (-5.234439323923323,-10.34747520197045)-- (-2.62044130297844,-11.012636345225395);
\draw [color=zzttqq] (-2.62044130297844,-11.012636345225395)-- (0.07678216708221663,-11.032898441709754);
\draw [color=zzttqq] (0.07678216708221663,-11.032898441709754)-- (2.700478262187136,-10.407083932025582);
\draw [color=zzttqq] (2.700478262187136,-10.407083932025582)-- (5.098167302036074,-9.171562880732537);
\draw [color=zzttqq] (5.098167302036074,-9.171562880732537)-- (7.130504314878069,-7.398139280395476);
\draw [color=zzttqq] (7.130504314878069,-7.398139280395476)-- (8.679377260771215,-5.189878064533485);
\draw [color=zzttqq] (8.679377260771215,-5.189878064533485)-- (9.654771274024288,-2.6751153483095353);
\draw [color=zzttqq] (9.7,-2.7)-- (10.,0.);
\draw [domain=-41.0:40.0] plot(\x,{(--1.68-0.2*\x)/22.37678575805597});
\draw [domain=-41.0:40.0] plot(\x,{(-28.64827273058606--1.2714034421419314*\x)/-2.3788565083480018});
\draw [domain=-41.0:40.0] plot(\x,{(--107.7-8.5*\x)/9.451344455077534});
\draw [domain=-41.0:40.0] plot(\x,{(--37.8-17.97332222129902*\x)/12.207768335121173});
\draw [domain=-41.0:40.0] plot(\x,{(-124.5-18.6*\x)/6.9});
\draw [domain=-41.0:40.0] plot(\x,{(--92.3--7.9*\x)/-0.9});
\draw [domain=-41.0:40.0] plot(\x,{(--92.5--7.9*\x)/1.0});
\draw [domain=-41.0:40.0] plot(\x,{(--125.6--18.5*\x)/7.2});
\draw [domain=-41.0:40.0] plot(\x,{(--36.0-17.8*\x)/-12.5});
\draw [domain=-41.0:40.0] plot(\x,{(-120.5--4.1*\x)/-16.2});
\draw [domain=-41.0:40.0] plot(\x,{(--106.3-8.4*\x)/-9.6});
\draw [domain=-41.0:40.0] plot(\x,{(-28.3--1.2*\x)/2.4});
\draw [domain=-41.0:40.0] plot(\x,{(-118.3--3.9*\x)/16.3});
\begin{scriptsize}
\draw (6.5,0.5) node[anchor=north west] {${Q_0}$};
\end{scriptsize}
\draw [fill=gray] (10.,0.) circle (10.5pt);
\draw [fill=ttqqqq] (9.7,2.7) circle (6.5pt);
\draw [fill=ttqqqq] (8.8,5.2) circle (6.5pt);
\draw [fill=ttqqqq] (7.2,7.4) circle (6.5pt);
\draw [fill=ttqqqq] (5.2,9.2) circle (6.5pt);
\draw [fill=ttqqqq] (2.92,10.5) circle (6.5pt);
\draw [fill=ttqqqq] (0.2,11.2) circle (6.5pt);
\draw [fill=ttqqqq] (-2.5,11.2) circle (6.5pt);
\draw [fill=ttqqqq] (-5.1,10.6) circle (6.5pt);
\draw [fill=ttqqqq] (-7.5,9.3) circle (6.5pt);
\draw [fill=ttqqqq] (-9.5,7.6) circle (6.5pt);
\draw [fill=ttqqqq] (-11.1,5.4) circle (6.5pt);
\draw [fill=ttqqqq] (-12.0,2.8) circle (6.5pt);
\draw [fill=ttqqqq] (-12.4,0.2) circle (6.5pt);
\draw [fill=ttqqqq] (-12.1,-2.5) circle (6.5pt);
\draw [fill=ttqqqq] (-11.13428252345037,-5.041033784216154) circle (6.5pt);
\draw [fill=ttqqqq] (-9.618760334548844,-7.272315365493406) circle (6.5pt);
\draw [fill=ttqqqq] (-7.613295832271324,-9.076071759828515) circle (6.5pt);
\draw [fill=ttqqqq] (-5.234439323923323,-10.34747520197045) circle (6.5pt);
\draw [fill=ttqqqq] (-2.62044130297844,-11.012636345225395) circle (6.5pt);
\draw [fill=ttqqqq] (0.07678216708221663,-11.032898441709754) circle (6.5pt);
\draw [fill=ttqqqq] (2.700478262187136,-10.407083932025582) circle (6.5pt);
\draw [fill=ttqqqq] (5.098167302036074,-9.171562880732537) circle (6.5pt);
\draw [fill=ttqqqq] (7.130504314878069,-7.398139280395476) circle (6.5pt);
\draw [fill=ttqqqq] (8.679377260771215,-5.189878064533485) circle (6.5pt);
\draw [fill=ttqqqq] (9.654771274024288,-2.6751153483095353) circle (12.5pt);
\draw [fill=ttqqqq] (-14.488379749344912,24.43103364331641) circle (12.5pt);
\draw [fill=ttqqqq] (-13.925551272742965,19.48552440628942) circle (12.5pt);
\draw [fill=ttqqqq] (-9.535753946027738,17.139353709610006) circle (12.5pt);
\draw [fill=ttqqqq] (-1.7618089573550433,12.984491748438458) circle (12.5pt);
\draw [fill=ttqqqq] (-6.749264649426915,0.1258239149020644) circle (12.5pt);
\draw [fill=ttqqqq] (-11.726557032080262,0.16321435718748709) circle (12.5pt);
\draw [fill=ttqqqq] (2.0650824138895443,0.059608730055132134) circle (12.5pt);
\draw [fill=ttqqqq] (-10.444793307539639,10.09537108236755) circle (12.5pt);
\draw [fill=ttqqqq] (-3.5617418313415765,8.34390135962629) circle (12.5pt);
\draw [fill=ttqqqq] (6.143498744449515,5.8742939363490585) circle (12.5pt);
\draw [fill=ttqqqq] (16.447574764747984,3.25230611121287) circle (12.5pt);
\draw [fill=ttqqqq] (12.69722347006066,-0.0202620964843597) circle (12.5pt);
\draw [fill=ttqqqq] (-14.852657156841106,-24.060373643680215) circle (12.5pt);
\draw [fill=ttqqqq] (22.71146313430742,-0.09549112090401835) circle (12.5pt);
\draw [fill=ttqqqq] (-14.215593032512771,-19.123878246225935) circle (12.5pt);
\draw [fill=ttqqqq] (-9.791043295911539,-16.84392273439346) circle (12.5pt);
\draw [fill=ttqqqq] (-1.9555548457597811,-12.806322042046812) circle (12.5pt);
\draw [fill=ttqqqq] (29.813659784221574,-0.14884428031732438) circle (12.5pt);
\draw [fill=ttqqqq] (16.39798586578574,-3.3488046000406255) circle (12.5pt);
\draw [fill=ttqqqq] (6.0556811656663045,-5.815692574217655) circle (12.5pt);
\draw [fill=ttqqqq] (-3.685566378030177,-8.139213986108459) circle (12.5pt);
\draw [fill=ttqqqq] (-10.594154347801913,-9.787078081717354) circle (12.5pt);
\end{tikzpicture}
\caption{. The $13$ lines of B\"or\"oczky.}
\label{b13conf}
\end{figure}

\section{Case $\mathbb{B}_{13}$}\label{Case13}

\subsection{Construction of $\mathbb{B}_{13}$}\label{ConstructionCase13}

The configurations $\mathbb{B}_{n}$ are the modifications of the  original construction of B\"or\"oczky on $n$ lines  preserving the  incidences of $\mathcal{B}_{n}$.

The original B\"or\"oczky construction of $\mathcal{B}_{13}$ uses trigonometric functions and is based on vertices of a regular $26$-gon. The core of the construction of $\mathbb{B}_{13}$
 is the set of four general points on the projective plane and five lines joining certain pairs of these
points. Since any two lines on the projective plane are projectively equivalent we may choose our points arbitrarily.
To simplify the calculations we take four standard points, i.e.,
$$P_1=(1:0:0),\hspace{0,5cm} P_2=(0:1:0),\hspace{0,5cm} P_3=(0:0:1),\hspace{0,5cm} P_4=(1:1:1).$$
Then we take the following lines
$$P_1P_4:\ z-y=0, \hspace{0,5cm}P_1P_2:\ z=0 ,\hspace{0,5cm} P_1P_3:\ y=0,$$ $$P_2P4:\ z-x=0, \hspace{0,5cm}P_3P_4:\ y-x=0.$$
It gives us the coordinates of intersection points
$$P_5=P_3P_4 \cap P_1P_2= (1:1:0)\ \textrm{ and }\ P_6= P_1P_3 \cap P_2P_4=(1:0:1).$$
From now on we need to introduce a parameter. We choose a point $P_7 \in P_1P_4$, distinct from all previous points.
We may  write its coordinates as
$$P_7= (a:1:1),$$
where  $a \neq 1$. Hence we obtain the following equations of lines
$$P_3P_7:\ ay-x=0, \hspace{1cm} P_2P_7:\ x-az=0.$$
In the next step we take the following points (and from now on we add the assumption that $a \neq 0$ which guarantees that the  points are distinct from the
already chosen ones).
$$
\begin{array}{ll}
P_8&=P_3P_7 \cap P_1P_2= (a:1:0),\\
P_9&=P_2P_7 \cap P_1P_3= (a:0:1),\\
P_{10}&=P_3P_4 \cap P_2P_7= (a:a:1),\\
P_{11}&=P_3P_7 \cap P_2P_4= (a:1:a).
\end{array}
$$
Now we choose an additional point $P_{12} \in P_1P_4$ distinct from all points from $P_1$ to $P_{11}$. Its coordinates are
$$P_{12}=(b:1:1),$$
where $b \neq 1$ and $b \neq a$. The next two lines of our construction are
$$P_{12}P_8:\ -x+ay+(b-a)z=0 \quad \textrm{  and } \quad P_{12}P_9:\ x+(a-b)y-az=0.$$
This gives us the following points
$$
\begin{array}{ll}
P_{13}&=P_{12}P_8 \cap P_2P_4= (a:a-b+1:a),\\
P_{14}&=P_{12}P_9 \cap P_3P_4= (a:a:a-b+1),\\
P_{15}&=P_1P_2 \cap P_{12}P_9= (b-a:1:0),\\
P_{16}&=P_{12}P_8 \cap P_1P_3= (b-a:0:1).
\end{array}
$$
The last four lines of the construction are
$$
\begin{array}{l}
P_{10}P_{15}:\ x-(b-a)y+(ab-a^2 -a)z=0,\\
P_{16}P_{11}:\ x+(ab-a^2 -a)y-(b-a)z=0,\\
P_{13}P_5:\ -ax+ay+(b-1)z=0,\\
P_6P_{14}:\ -ax+(b-1)y+az=0.
\end{array}
$$
Finally we obtain the remaining triple points
$$
\begin{array}{lll}
P_{17}&=P_1P_4 \cap P_{10}P_{15} \cap P_{11}P= (a^2+b-ab:1:1),\\
P_{18}&=P_{12}P_8 \cap P_{10}P_{15} \cap P_6P_{14}\\&= (b^2-a^2-b+a-ab:ab-a^2-a:-a^2+b-1),\\
P_{19}&=P_2P_7 \cap P_6P_{14} \cap P_{11}P= (a-ab:a-a^2:1-b),\\
P_{20}&=P_{12}P_9 \cap P_{13}P_5 \cap P_{11}P\\&= (b^2-a^2-b+a-ab:-a^2+b-1:ab-a^2-a),\\
P_{21}&=P_{13}P_5 \cap P_{10}P_{15} \cap P_3P_7= (a-ab:a-a^2:1-b),\\
P_{22}&=P_1P_4 \cap P_{13}P_5 \cap P_6P_{14}= (a+b-1:a:a).
\end{array}
$$

It is easy to check (computing a suitable determinant) that points $P_{17}$ and $P_{22}$ are always the common points of given three lines, independently of
the values of $a$ and $b$. Situation is more complicated for points $P_{18}$, $P_{19}$, $P_{20}$ and $P_{21}$. These points are triple only if parameters $a$
and $b$ satisfy the additional condition:
$$C_{13}(a,b):  a^4-a^3b+a^2b -a^2+b^2-2ab+2a-b=0.$$
This condition is necessary for the construction to terminate successfully, in the sense that we obtain exactly $22$ triple points on $13$ lines,
lying in groups: of six points on  one  of the lines and of five points on each of the remaining 12 lines, as in the original B\"or\"oczky configuration $\mathcal{B}_{13}$.

\subsection{Degenerate cases of $\mathbb{B}_{13}$}\label{DegenerateCase13}

In this section we check under what conditions the configuration $\mathbb{B}_{13}$ in non-degene\-ra\-te, in the sense that all points and lines appearing in
the construction are distinct. Due to the choice of points $P_{7}$ and $P_{12}$ we are working with the assumptions $a \neq 0$, $a \neq 1$, $a \neq b$ and $b
\neq 1$. These conditions force that the chosen points $P_{7}$ and $P_{12}$ are distinct from the four points starting the construction and $P_{7} \neq
P_{12}$. There is also the prevail  condition over any others, guaranteeing  that construction contains exactly $22$ triple points, namely
\begin{equation} \label{main13}
C_{13}(a,b):  a^4-a^3b+a^2b -a^2+b^2-2ab+2a-b=0.
\end{equation}

In spite of the above  conditions, by comparing the equations of lines and coordinates of all our points, we obtain new terms, saying when some of points and
lines overlap.

\begin{itemize}
\item[$i)$] if $a-b+1=0$ then $P_{4}=P_{17}$, $P_{5}=P_{14}=P_{15}=P_{20}=P_{10}$ and $P_{6}=P_{13}=P_{18}=P_{16}=P_{11}$,
\item[$ii)$] if $2a-b=0$ then $P_{8}=0$, $P_{9}=P_{16}$, $P_{11}=P_{20}=P_{13}$ and $P_{12}=P_{18}=P_{14}$,
\item[$iii)$] if $a^2-ab+2a-b=0$ then $P_{10}P_{15}=P_{16}P_{11}$ (thus $P_{10},P_{15},P_{16},P_{11}$ are collinear points),
\item[$iv)$] if $a^2-b+1=0$ then $P_{8}=P_{18}$ and $P_{9}=P_{20}$,
\item[$v)$] if $a^2-a-b+1=0$ then  $P_{7}=P_{14}=P_{21}=P_{22}=P_{19}$,
\item[$vi)$] if $a^3-a^2b+ab-a-b+1=0$ then $P_{17}=P_{22}$,
\item[$vii)$] if $a^3-2ab+2a+b^2-a^2-b=0$ then $P_{21}=P_{20}$ and $P_{18}=P_{19}$.
\end{itemize}

Together with main condition (\ref{main13}) it gives us
\begin{itemize}
\item[$i)$] $-a(a-1)=0$
\item[$ii)$] $-a^2(a-1)^2=0$
\item[$iii)$] $- \frac{a^3}{(a+1)^2}=0$
\item[$iv)$] $-a^2(a-1)^3=0$
\item[$v)$] $-a (a-1)^4=0$
\item[$vi)$] $ \frac{a^2 (a-1)^3}{(a^2-a+1)^2}=0$
\item[$vii)$] $a^2(a-1)(a-b)=0$
\end{itemize}
It is easy to see, that all cases give us solutions which contradict with permitted values of coefficients $a$ and $b$. Thus, if only  $a \neq 0$, $a \neq 1$,
$a \neq b$ and $b \neq 1$, then the construction always consists of $13$ distinct lines, which give the $22$ triple intersection points (provided that
$C_{13}(a,b)=0)$.

\subsection{Parameter space of configurations $\mathbb{B}_{13}$}\label{ParameterSpaceCase13}

The parametrizing curve
$$C_{13}(a,b): a^4-a^3b+a^2b -a^2+b^2-2ab+2a-b=0$$
is an irreducible curve of degree $4$ with one double point. This may be checked by hand or with the help of a computer. Thus, the geometrical genus of $C_{13}$ is
$2$. The curve may be written as
$$ b^2+b(-1-2a+a^2-a^3)+a^4-a^2+2a=0.$$
Substituting $b$ by
$b-\frac{-1-2a+a^2-a^3}{2}$ and then $b$ by $\frac{b}{2}$ we get the curve
$$D_{13}: \ 1 - 4a + 6a^2 - 2a^3 + a^4 - 2a^5 + a^6 -b^2=0.$$
We will denote the homogenization of $D_{13}$ also by $D_{13}$.
Using computer algebra programme Magma \cite{Magmaa}, it may be computed  that the Jacobian of $D_{13}$ has rank $0$, thus Chabauty's method may be applied here (cf. e.g. [Poonen]).
We obtain all rational points of $D_{13}$, namely
$$(1 : -1 : 0), \hspace{0,3cm}(1 : 1 : 0), \hspace{0,3cm}(0 : 1 : 1), \hspace{0,3cm}(0 : -1 : 1), \hspace{0,3cm}(1 : 1 : 1), \hspace{0,3cm}(1 : -1 : 1).$$

These points correspond to the set of all rational points on $C_{13}$ (with perhaps second coordinate changed). However,  points having $a=1$ or $a=0$ are excluded
by the construction. Thus all possibilities lead us to the degenerated cases. Thus we have the following

\begin{corollary}
The configuration  $\mathbb{B}_{13}$ cannot be realized over the rational numbers.
\end{corollary}

\section{Construction of $\mathbb{B}_{14}$, $\mathbb{B}_{16}$, $\mathbb{B}_{18}$ and $\mathbb{B}_{24}$}\label{14-24}

For  configurations  $\mathbb{B}_{14}$, $\mathbb{B}_{16}$, $\mathbb{B}_{18}$ and $\mathbb{B}_{24}$ we present simplified descriptions of construction in
tables.  We omit the coordinates of points and the equations of lines of the configurations, because of their complicated forms. We only give the coordinates
of points being the core of each construction and coordinates of points taken during the construction as  parameters. We distinguish these special points using bold type font. Enough motivated reader may reconstruct the construction step by step and find remaining coordinates and equations of lines if  necessary. 

The idea of the construction is the same in all considered cases: we start with four fundamental points and some of the lines through them, and then we choose two
points (pa\-ra\-meters) on one of the already constructed lines. This input allows to construct the configurations.

\subsection{Case of $ \mathbb{B}_{14}$}\label{Case14}

The construction of $\mathbb{B}_{14} $ is based on four standard points on the projective plane. We introduce here two parameters $a$ and $b$ such that $a\neq
1$, $b\neq 1$ and $a \neq b$. The construction runs in the following way
\begin{center}
\begin{tabular}{c|l}
    step 1 & $P_1=(1:0:0)$, $B=(0:1:0)$, $P_3=(0:0:1)$, $P_4=(1:1:1)$\\
    \hline
    step 2 &  lines: $P_1P_2$, $P_1P_3$, $P_1P_4$, $P_2P_3$, $P_2P_4$, $P_3P_4$\\
    \hline
    step 3& $P_5=P_1P_2 \cap P_3P_4$, $P_6= P_1P_3 \cap P_2P_4$, ${\bf P_7= (a:1:1) \in P_1P_4}$\\
    \hline
    step 4  &  lines: $P_5P_7$, $P_6P_7$ \\
    \hline
    \multirow{3}{1cm}{step 5}&$P_8=P_2P_4 \cap P_5P_7$, $P_9=P_3P_4 \cap P_6P_7$, $P_{10}=P_1 P_2\cap P_6P_7$,\\
    & $P_{11}=P_1P_3 \cap P_5P_7$, $P_{12}=P_2P_3 \cap P_5P_7$, $P_{13}=P_2P_3 \cap P_6P_7$,\\
    & ${\bf P_{14}=(b:1:1) \in P_1P_4}$\\
    \hline
    step 6 &lines: $P_{10}P_{14}$, $P_{11}P_{14}$\\
    \hline
    \multirow{3}{1cm}{step 7} &  $P_{15}=P_5P_7 \cap P_{10}P_{14},$
    $P_{16}=P_3P_4 \cap P_{10}P_{14},$
    $P_{17}=P_2P_4 \cap P_{10}P_{14},$\\
    &$P_{18}=P_2P_3 \cap P_{10}P_{14},$
     $P_{19}=P_6P_7 \cap P_{11}P_{14},$
    $P_{20}=P_{2}P_{4} \cap P_{11}P_{14},$\\
    &$P_{21}=P_3P_4 \cap P_{11}P_{14},$
    $P_{22}=P_2P_3 \cap P_{11}P_{14}$ \\
    \hline
        step 8 &lines: $P_{15}P_{20}$, $P_{16}P_{19}$, $P_{17}P_{22}$, $P_{18}P_{21}$\\
    \hline
        \multirow{2}{1cm}{step 9} & $P_{23}=P_1P_4 \cap P_{15}P_{20} $,
        $P_{24}=P_1P_3 \cap P_{15}P_{20} $,
        $P_{25}=P_1P_2 \cap P_{16}P_{19}$,\\
        & $P_{26}=P_{17}P_{22} \cap P_{18}P_{21} $\\
\end{tabular}

\end{center}

There exist some points, which are not always triple, namely $P_9$, $P_{12}$, $P_{13}$, $P_{24}$ and $P_{26}$. The conditions to have them triple  (i.e., to
assure additional incidences) are
$$
\begin{array}{ll}
P_9 &\in P_3P_4 \cap P_6P_7 \cap P_{17}P_{22},\\
P_{12} &\in P_2P_3 \cap P_5P_7 \cap P_{16}P_{19},\\
P_{13} &\in P_2P_3 \cap P_6P_7 \cap P_{15}P_{20},\\
P_{24} &\in P_1P_3 \cap P_{18}P_{21} \cap P_{15}P_{20},\\
P_{26} &\in P_1P_4 \cap P_{18}P_{21} \cap P_{17}P_{22},
\end{array}
$$
what implies the following algebraic conditions for parameters $a$ and $b$:

\begin{itemize}
\item{}$2 a - 3 a^2 + a^3 + a b - b^2=0$,
\item{} $(-2 + a) (2 a - 3 a^2 + a^3 + a b - b^2)=0$,
\item{}$(-2 + 2 a - b) (-1 + b) (2 a - 3 a^2 + a^3 + a b - b^2)=0$.
\end{itemize}
Thus the  curve parametrizing the configurations $\mathbb{B}_{14}$ is the curve
$$C_{14}(a,b): 2 a - 3 a^2 + a^3 + a b - b^2=0.$$

\subsubsection{Parameter space of configurations $\mathbb{B}_{14}$} We now take a closer look at the curve parametrizing the configurations $\mathbb{B}_{14}$,
i.e.,
$$C_{14}(a,b):  2 a - 3 a^2 + a^3 + a b - b^2=0.$$
Its geometric genus is $1$.  By substituting $a\mapsto a/4, b\mapsto (a+b)/4$ we get
$$D_{14}: b^2=a^3-11a^2+32a=0.$$
The curve $D_{14}$ has the following rational points:
$$(0,0),\hspace{0,5cm}(4,-4),\hspace{0,5cm}(4,4),\hspace{0,5cm}(8,-8),\hspace{0,5cm}(8,8).$$
Thus, $C_{14}$ has on the projective plane the following rational points
$$(0:0:1),  \hspace{0,5cm}(1:1:1), \hspace{0,5cm}(1:0:1), \hspace{0,5cm}(2:2:1), \hspace{0,5cm}(2:0:1), \hspace{0,5cm}(0:1:0).$$
 Each of them leads us to a degenerated case hence no $\mathbb{B}_{14}$ configuration  can be obtained over the
 rational numbers.

\subsection{Case of $\mathbb{B}_{16}$}
The core of configuration $\mathbb{B}_{16}$ are four general points. We need two parameters
$a$ and $b$ such that $a\neq 1$, $a\neq 0$, $b\neq 0$, $b\neq1$ and $a\neq b$. We present the construction step by step.
\begin{center}
    \begin{tabular}{c|l}
        step 1 &
        $P_1=(1:0:0)$, $P_2=(0:1:0)$, $P_3=(0:0:1)$, $P_4=(1:1:1)$\\
        \hline
        step 2 &  lines: $ P_1 P_4 $,
        $ P_2 P_4 $,
        $ P_3 P_4 $ \\
        \hline
    {step 3}&
        ${\bf P_5 =(a:1:1) \in P_1 P_4}$\\
        \hline
        step 4  &  lines: $P_2 P_5$, $P_3 P_5$ \\
        \hline
    {step 5}&$P_{6} =P_2P_4 \cap P_3P_5  $,
        $P_{7} =P_3 P_4 \cap P_2 P_5  $\\
        \hline
        step 6 &lines: $P_1 P_{6}$, $P_1 P_{7}$\\
        \hline
        \multirow{2}{1cm}{step 7} &  $P_{8} =P_3 P_5\cap P_1 P_{7} , $
        $P_{9} =P_2 P_4  \cap  P_1 P_{7} , $
        $P_{10} =P_1 P_6\cap P_2 P_{5} , $\\
        & $P_{11} =P_3 P_4\cap P_1 P_{6}  $\\
        \hline
        step 8 &line $P_{9} P_{11}$\\
        \hline
    {step 9} & $P_{12} =P_9 P_{11} \cap P_{1} P_{5}  , $
        $P_{13} =P_9 P_{11} \cap P_{2} P_{5}, $ ${\bf P_{14} =(b:1:1) \in P_1 P_4}$ \\
        \hline
            step 10 &lines: $P_{2} P_{14}$, $P_{3} P_{14}$\\
        \hline
        \multirow{3}{1,2cm}{step 11} & $P_{15} =P_3 P_{14} \cap P_{1} P_{6} $,
        $P_{16} =P_3 P_{14} \cap P_{9} P_{11}$,
        $P_{17} =P_3 P_{14}  \cap P_{1} P_{7} $,\\
        & $P_{18} =P_2 P_{14} \cap P_{3} P_{4} $,
        $P_{19} = P_2 P_{14}\cap P_{1} P_{6}$, $P_{20} = P_2 P_{14}\cap P_{9} P_{11} $\\
        &$P_{21} = P_2 P_{14}\cap P_{1} P_{7} $, $P_{22} =P_{2} P_{4} \cap P_{3} P_{14} $\\
        \hline
            step 12 &lines: $P_8P_{20}$, $P_{10}P_{16}$, $P_{12} P_{19}$,  $P_{13}P_{17}$,  $P_{15} P_{18}$,   $P_{21} P_{22}$\\
    \hline
\multirow{5}{1,2cm}{step 13} & $P_{23} =P_3 P_{14} \cap P_{8} P_{20} $,
$P_{24} =P_{15} P_{18} \cap P_{13} P_{17} ,$
$P_{25} =P_{13} P_{17}  \cap P_{8} P_{20} $,\\
& $P_{26} =P_{13} P_{17} \cap P_{1} P_{4} $, $P_{27}=P_2P_4 \cap P_{15}P_{18}$,
    $P_{28}=P_{10}P_{16} \cap P_{13}P_{17}$,\\
&  $P_{29}=P_2P_4 \cap P_{10}P_{16}$,
        $P_{30}=P_{2}P_{5} \cap P_{12}P_{19}$, $P_{31} =P_{10} P_{16} \cap P_{3} P_{4} $,\\
&  $P_{32}=P_1P_4 \cap P_{15}P_{18}$,
        $P_{33}=P_{1}P_{6} \cap P_{8}P_{20}$, $P_{34} =P_{1} P_{7} \cap P_{10} P_{16} $,\\
&  $P_{35}=P_1P_4 \cap P_{10}P_{16}$
    \end{tabular}

\end{center}
There exist some points, which are not always triple, namely $P_{23}$, $P_{24}$, $P_{25}$, $P_{28}$, $P_{30}$, $P_{33}$ and $P_{34}$.
The conditions to have them triple are
$$
\begin{array}{ll}
P_{23}&=P_{3}P_{14} \cap P_{8}P_{20} \cap P_{12}P_{19},\\
P_{24}&=P_{3}P_{6} \cap P_{13}P_{17} \cap P_{15}P_{18},\\
P_{25}&=P_{2}P_{4} \cap P_{8}P_{20} \cap P_{15}P_{18},\\
P_{28}&=P_{2}P_{14} \cap P_{10}P_{16} \cap P_{13}P_{17},\\
P_{30}&=P_{2}P_{5} \cap P_{12}P_{19} \cap P_{21}P_{22},\\
P_{33}&= P_{1}P_{6} \cap P_{8}P_{20} \cap P_{21}P_{22},\\
P_{34}&=P_{1}P_{7} \cap P_{10}P_{16} \cap P_{15}P_{18}.
\end{array}
$$

It implies the condition for $a$ and $b$
$$C_{16}(a,b):=
a^4b^2-2a^3b-2a^2b^2+a^2b+2ab^2-b^3+a^2+2b^2-2b.$$

\subsubsection{Parameter space of configurations $\mathbb{B}_{16}$}

The parametrizing curve is
$$C_{16}(a,b):\ 
a^4b^2-2a^3b-2a^2b^2+a^2b+2ab^2-b^3+a^2+2b^2-2b=0.$$ Computations with Magma returned genus of $C_{16}=2$. This curve is hyperelliptic and its Jacobian 
$J$ has rank $0$ (RankBound($J$)$=0$, as computed by Magma). Thus, Chabauty's method may be applied  and Magma computes that the  rational points  on $C_{16}$
 are:
 $$
  (1 : 1 : 1), \hspace{10pt} (0 : 1 : 0),\hspace{10pt} (0 : 0 : 1), \hspace{10pt} (1 : 0 : 0),$$ 
$$(-2 : -2 : 1),\hspace{10pt} (-1 : 1 : 1),\hspace{10pt} (-1 : -1 : 1).
$$
All these points give  degenerated configurations.

\subsection{ Case of $\mathbb{B}_{18}$}
The construction of $\mathbb{B}_{18}$ is also based on four standard points on the projective plane. It also depends of parameters $a$ and $b$ satisfying the conditions
$a\neq 1$, $b\neq 1$ and $a\neq b$.

\begin{center}
    \begin{tabular}{c|l}

        step 1 &
        $P_1=(1:0:0)$, $P_2=(0:1:0)$, $P_3=(0:0:1)$, $P_4=(1:1:1)$\\
        \hline
        step 2 &  lines: $ P_1 P_4 $,
        $ P_2 P_4 $,
        $ P_3 P_4 $ \\
        \hline
        {step 3}&
        ${\bf P_5 =(a:1:1) \in P_1 P_4}$\\
        \hline
        step 4  &  lines: $P_2 P_5$, $P_3 P_5$ \\
        \hline
        {step 5}&$P_{6} =P_2P_4 \cap P_3P_5  $,
        $P_{7} =P_3 P_4 \cap P_2 P_5  $\\
        \hline
        {step 6}&
        ${\bf P_8 =(b:1:1) \in P_1 P_4}$\\
        \hline
        step 7  &  lines: $P_6 P_8$, $P_7 P_8$ \\
        \hline
        \multirow{2}{1,1cm}{step 8} &  $P_{9} =P_2 P_5\cap P_6 P_{8} , $
        $P_{10} =P_3 P_5  \cap  P_7 P_{8} , $
        $P_{11} =P_3 P_4\cap P_6 P_{8} , $\\
        & $P_{12} =P_2 P_4\cap P_7 P_{8}  $\\
        \hline
        step 9 &line $P_{1} P_{9}$,  $P_{1} P_{10}$\\
        \hline
        \multirow{2}{1,3cm}{step 10} &  $P_{13} =P_3 P_4\cap P_1 P_{9} , $
        $P_{14} =P_2 P_4  \cap  P_1 P_{10} , $
        $P_{15} =P_1 P_{10}\cap P_6 P_{8} , $\\
        & $P_{16} =P_1 P_9\cap P_7 P_{8}  $, $P_{17}=P_2P_4 \cap P_1P_9$,
        $P_{18}=P_3P_4 \cap P_1P_{10}$\\
        \hline
        step 11 &lines: $P_{13} P_{15}$, $P_{14} P_{16}$\\
        \hline
        \multirow{3}{1,3cm}{step 12} & $P_{19} =P_2 P_{4} \cap P_{13} P_{15}$,
        $P_{20} =P_3 P_{4} \cap P_{14} P_{15} ,$
        $P_{21} =P_2 P_{5}  \cap P_{13} P_{15} $,\\
        & $P_{22} =P_3 P_{5} \cap P_{14} P_{16} $,
        $P_{23} = P_1 P_{4}\cap P_{13} P_{15} $, $P_{24} = P_7 P_{10}\cap P_{13} P_{15} $\\
        &$P_{25} = P_6 P_{9}\cap P_{14} P_{16} $\\
        \hline
        step 13 &lines: $P_{17} P_{21}$, $P_{18} P_{22}$, $P_{19}P_{20}$\\
        \hline
        \multirow{5}{1,3cm}{step 14} & $P_{26} =P_{17} P_{21} \cap P_{1} P_{4} $,
        $P_{27} =P_{3} P_{5} \cap P_{17} P_{21} ,$
        $P_{28} =P_{18} P_{22}  \cap P_{2} P_{5} $,\\
        & $P_{29} =P_{19} P_{20} \cap P_{17} P_{21} $, $P_{30}=P_{19}P_{20} \cap P_{18}P_{22}$,
        $P_{31}=P_{3}P_{4} \cap P_{17}P_{21}$,\\
        & $P_{32} =P_{2} P_{4} \cap P_{18} P_{22} $, $P_{33}=P_1P_{10} \cap P_{17}P_{21}$,
        $P_{34}=P_{1}P_{9} \cap P_{18}P_{22}$,\\
        & $P_{35} =P_{17} P_{21} \cap P_{14} P_{16} $, $P_{36}=P_{18}P_{22} \cap P_{13}P_{15}$,
        $P_{37}=P_{19}P_{20} \cap P_{1}P_{9}$,\\
        & $P_{38} =P_{19} P_{20} \cap P_{1} P_{10} $, $P_{39}=P_3P_5 \cap P_{11}P_{20}$,
        $P_{40}=P_2P_5 \cap P_{19}P_{20}$\\
        \hline
        step 15 &lines: $P_{2} P_{24}$, $P_{3} P_{25}$, $P_{11} P_{27}$, $P_{12} P_{28}$    \\
        \hline
        \multirow{2}{1,3cm}{step 16} & $P_{41} =P_{2} P_{24} \cap P_{1} P_{10} $,
        $P_{42} =P_{3} P_{25} \cap P_{1} P_{9} ,$
        $P_{43} =P_{1} P_{4}  \cap P_{11} P_{27} $,\\
        & $P_{44} =P_{2} P_{24} \cap P_{6} P_{8} $, $P_{45}=P_{3}P_{25} \cap P_{7}P_{8}$,
        $P_{46}=P_{1}P_{4} \cap P_{2}P_{24}$
    \end{tabular}
\end{center}
\normalsize

We obtain $46$ triple intersection points. Here  there also exist some points which are not necessarily triple. The conditions to have them triple are:

\begin{minipage}{0.45\linewidth}
$$
\begin{array}{ll}
P_{29}&=P_{6}P_{8} \cap P_{17}P_{21} \cap P_{19}P_{20},\\
P_{30}&=P_{7}P_{8} \cap P_{18}P_{22} \cap P_{19}P_{20},\\
P_{31}&=P_{2}P_{24} \cap P_{3}P_{4} \cap P_{17}P_{21},\\
P_{32}&=P_{2}P_{4} \cap P_{3}P_{25} \cap  P_{18}P_{22},\\
P_{35}&= P_{11}P_{27}\cap P_{12}P_{28}\cap P_{14}P_{16} ,\\
P_{36}&= P_{11}P_{27}\cap P_{13}P_{15} \cap P_{18}P_{22},
\end{array}
$$
\end{minipage}
\hskip 20pt
\begin{minipage}{0.4\linewidth}
$$
\begin{array}{ll}
P_{37}&=P_{1}P_{9} \cap P_{11}P_{27} \cap P_{19}P_{20},\\
P_{38}&=P_{1}P_{10} \cap P_{12}P_{28} \cap P_{19}P_{20},\\
P_{39}&=P_{2}P_{24} \cap P_{3}P_{5} \cap P_{19}P_{20},\\
P_{40}&=P_{2}P_{5} \cap P_{3}P_{25} \cap  P_{19}P_{20},\\
P_{44}&=P_{2}P_{24} \cap P_{6}P_{8} \cap P_{12}P_{28},\\
P_{45}&=P_{3}P_{25} \cap P_{7}P_{8} \cap P_{11}P_{27}.
\end{array}
$$
\end{minipage}

\vskip 5pt
They imply several algebraic conditions for parameters $a$ and $b$ but only one of them does not lead to a degenerated case. This condition provides us the parametrization curve of the $\mathbb{B}_{18}$ configurations
\begin{equation*}
\begin{aligned}
C_{18}(a,b):=& a^3b^5 -a^5b^2 +a^4b^3-6a^3b^4 + a^2b^5+a^6+a^4b^2+12a^3b^3-4a^2b^4-\\&5a^5+ 7a^4b-22a^3b^2+11a^2b^3-ab^4 + 6a^4-a^3b+3a^2b^2 -4ab^3+ \\&b^4-4a^3+ 4a^2b-ab^2.
\end{aligned}
\end{equation*}

\subsubsection{Parameter space of configurations $\mathbb{B}_{18}$}

Computations with Magma returned  genus of 
\begin{equation*}
\begin{aligned}
C_{18}(a,b):=& a^3b^5 -a^5b^2 +a^4b^3-6a^3b^4 + a^2b^5+a^6+a^4b^2+12a^3b^3-4a^2b^4-\\&5a^5+ 7a^4b-22a^3b^2+11a^2b^3-ab^4 + 6a^4-a^3b+3a^2b^2 -4ab^3+ \\&b^4-4a^3+ 4a^2b-ab^2
\end{aligned}
\end{equation*}
 is equal to $2$ again. This curve is hyperelliptic and its Jacobian 
$J$ has again rank $0$. Thus, Chabauty's method may be applied  and Magma computes the  rational points  on $C_{18}$:
$$
 (1 : 1 : 1),\hspace{10pt}  (0 : 1 : 0),\hspace{10pt} (0 : 0 : 1),\hspace{10pt} (1 : 0 : 0).
$$
All these points give a degenerated configuration.

\subsection{Construction of $\mathbb{B}_{24}$}

We start the construction with taking four standard points. Then we introduce  the parameters $a$ and $b$, such that
$a\neq 1$, $a\neq 0$, $b\neq1$ and $a\neq b$ and we have
\begin{center}
    \begin{tabular}{c|l}
        step 1 &
        $P_1=(1:0:0)$, $P_2=(0:1:0)$, $P_3=(0:0:1)$, $P_4=(1:1:1)$\\
        \hline
        step 2 &  lines: $P_1P_4$, $P_2P_4$, $ P_3P_4$\\
        \hline
    {step 3}& ${\bf P_5 =(a:1:1) \in P_1 P_4}$\\
        \hline
        step 4  &  lines: $P_2P_5$, $P_3P_5$ \\
        \hline
    {step 5}& $P_6= P_2P_4 \cap  P_3P_5,$ $P_7= P_2P_5 \cap  P_3P_4$\\
        \hline
        step 6 &lines: $P_1P_6$, $P_1P_7$\\
        \hline
        {step 7} &  $P_8= P_1P_6 \cap  P_3P_4,$
        $P_9= P_1P_7 \cap  P_2P_4$  \\
        \hline
        step 8 &lines: $P_2P_8$, $P_3P_9$\\
        \hline
        \multirow{3}{1cm}{step 9} & $P_{10}= P_1P_4 \cap  P_3P_9,$
        $P_{11}= P_2P_5 \cap  P_3P_9,$
        $P_{12}= P_2P_8 \cap  P_3P_5,$\\
        &$P_{13}= P_1P_6 \cap  P_2P_5,$
        $P_{14}= P_1P_7 \cap  P_3P_5,$
        $P_{15}= P_1P_6 \cap  P_3P_9,$\\
        &$P_{16}= P_1P_7 \cap  P_2P_8$\\
        \hline
        step 10 & ${\bf P_{17}=(b:1:1)\in P_1P_4}$\\
        \hline
        step 11 & lines: $P_{11}P_{17}$, $P_{12}P_{17}$\\
        \hline
        \multirow{4}{1,2cm}{step 12} & $P_{18}= P_3P_9 \cap  P_{12}P_{17},$
        $P_{19}= P_2P_8 \cap  P_{11}P_{17},$
        $P_{20}= P_1P_7 \cap  P_{11}P_{17},$\\
        &$P_{21}= P_1P_6 \cap  P_{12}P_{17},$
        $P_{22}= P_3P_4 \cap  P_{12}P_{17},$
        $P_{23}= P_2P_4 \cap  P_{11}P_{17},$\\
        &$P_{24}= P_2P_5 \cap  P_{12}P_{17},$
        $P_{25}= P_3P_5 \cap  P_{11}P_{17},$
        $P_{26}= P_3P_4 \cap  P_{11}P_{17},$\\
        &$P_{27}= P_2P_4 \cap  P_{12}P_{17}$\\
        \hline
        step 13 &lines: $P_{12}P_{22} $, $P_{22}P_{23}$\\
        \hline
        \multirow{4}{1,2cm}{step 14} & $P_{28}= P_4P_9 \cap  P_{22}P_{23},$
        $P_{29}= P_{12}P_{17} \cap  P_{22}P_{23},$
        $P_{30}= P_{2}P_{5} \cap  P_{22}P_{23},$\\
        & $P_{31}= P_{3}P_{4} \cap  P_{22}P_{23},$
        $P_{32}= P_{1}P_{45} \cap  P_{22}P_{23},$
        $P_{33}= P_{1}P_{6} \cap  P_{22}P_{23},$\\
        & $P_{34}= P_{1}P_{7} \cap  P_{15}P_{22},$
        $P_{35}= P_{2}P_{4} \cap  P_{15}P_{22},$
        $P_{36}= P_{3}P_{5} \cap  P_{15}P_{22},$\\
        & $P_{37}= P_{11}P_{17} \cap  P_{15}P_{22},$
        $P_{38}= P_{2}P_{8} \cap  P_{15}P_{22}$\\
        \hline
            step 15 &lines: $P_{24}P_{28}  $, $P_{25}P_{38}$\\
        \hline
        \multirow{5}{1,2cm}{step 16} & $P_{39}= P_{1}P_{4} \cap  P_{24}P_{28},$
        $P_{40}= P_{2}P_{4} \cap  P_{24}P_{28},$
        $P_{41}= P_{1}P_{7} \cap  P_{24}P_{28},$\\
        & $P_{42}= P_{2}P_{8} \cap  P_{24}P_{28},$
        $P_{43}= P_{3}P_{4} \cap  P_{24}P_{28},$
        $P_{44}= P_{15}P_{22} \cap  P_{24}P_{28},$\\
        & $P_{45}= P_{1}P_{6} \cap  P_{24}P_{28},$
        $P_{46}= P_{3}P_{5} \cap  P_{24}P_{28},$
        $P_{47}= P_{2}P_{5} \cap  P_{25}P_{38},$\\
        & $P_{48}= P_{1}P_{7} \cap  P_{25}P_{38},$
        $P_{49}=P_{22}P_{23} \cap  P_{25}P_{38},$
        $P_{50}= P_{2}P_{4} \cap  P_{25}P_{38},$\\
        & $P_{51}= P_{3}P_{9} \cap  P_{25}P_{38},$
        $P_{52}= P_{1}P_{6} \cap  P_{25}P_{38},$
        $P_{53}= P_{3}P_{4} \cap  P_{25}P_{38}$\\
        \hline
            step 17 &lines: $P_{30}P_{53}  $, $P_{36}P_{40}$\\
        \hline
        \multirow{5}{1,2cm}{step 18} & $P_{54}= P_{3}P_{9} \cap  P_{36}P_{40},$
        $P_{55}= P_{2}P_{5} \cap  P_{36}P_{40},$
        $P_{56}= P_{2}P_{8} \cap  P_{36}P_{40},$\\
        &$P_{57}= P_{12}P_{17} \cap  P_{36}P_{40},$
        $P_{58}= P_{1}P_{4} \cap  P_{36}P_{40},$
        $P_{59}= P_{3}P_{4} \cap  P_{36}P_{40},$\\
        &$P_{60}= P_{25}P_{38} \cap  P_{36}P_{40},$
        $P_{61}= P_{1}P_{6} \cap  P_{36}P_{40},$
        $P_{62}= P_{1}P_{7} \cap  P_{30}P_{53},$\\
        &$P_{63}= P_{24}P_{28} \cap  P_{30}P_{53},$
        $P_{64}= P_{2}P_{4} \cap  P_{30}P_{53},$
        $P_{65}= P_{11}P_{17} \cap  P_{30}P_{53},$\\
        &$P_{66}= P_{3}P_{9} \cap  P_{30}P_{53},$
        $P_{67}= P_{3}P_{5} \cap  P_{30}P_{53},$
        $P_{68}= P_{2}P_{8} \cap  P_{30}P_{53}$\\
        \hline
            step 19 &lines: $P_{42}P_{54}  $, $P_{51}P_{68}$\\
        \hline
        \multirow{3}{1,2cm}{step 20} & $P_{69}= P_{2}P_{4} \cap  P_{42}P_{54},$
        $P_{70}= P_{15}P_{22} \cap  P_{42}P_{54},$
        $P_{71}= P_{1}P_{6} \cap  P_{42}P_{54},$\\
        &$P_{72}= P_{1}P_{4} \cap  P_{42}P_{54},$
        $P_{73}= P_{3}P_{5} \cap  P_{42}P_{54},$
        $P_{74}= P_{2}P_{5} \cap  P_{51}P_{68},$\\
        &$P_{75}= P_{1}P_{7} \cap  P_{51}P_{68},$
        $P_{76}= P_{22}P_{23} \cap  P_{51}P_{68},$
        $P_{77}= P_{3}P_{4} \cap  P_{51}P_{68}$\\
        \hline
            step 21 &lines: $P_{26}P_{66}  $, $P_{27}P_{56}$\\
        \hline
        \multirow{2}{1,2cm}{step 22} & $P_{78}= P_{11}P_{17} \cap  P_{27}P_{56},$
        $P_{79}= P_{3}P_{5} \cap  P_{27}P_{56},$
        $P_{80}= P_{1}P_{4} \cap  P_{27}P_{56},$\\
        &$P_{81}= P_{2}P_{5} \cap  P_{26}P_{66},$
        $P_{82}= P_{12}P_{17} \cap  P_{26}P_{66}$\\
        \hline
            step 23 &lines: $P_{18}P_{65} $, $P_{19}P_{57}$\\
        \hline
        {step 24} & $P_{83}= P_{3}P_{9} \cap  P_{19}P_{57},$
        $P_{84}= P_{2}P_{8} \cap  P_{18}P_{65},$
        $P_{85}= P_{18}P_{65} \cap P_{19}P_{57}$\\
        \hline
        step 25 &line: $P_{41}P_{52}$\\
    \end{tabular}
\end{center}
\normalsize

In this case we obtain $85$ triple intersection points and, as previously, some points are not necessarily triple. The conditions to have them triple are:

\begin{minipage}{0.45\linewidth}
$$
\begin{array}{ll}
P_{29}&=P_{12}P_{17} \cap  P_{22}P_{23}\cap P_{42}P_{54},\\
P_{33}&= P_{1}P_{6} \cap  P_{18}P_{65} \cap  P_{22}P_{23},\\
P_{34}&=P_{1}P_{7} \cap  P_{15}P_{22} \cap  P_{19}P_{57},\\
P_{37}&= P_{11}P_{17} \cap  P_{15}P_{22}\cap P_{51}P_{68},\\
P_{43}&=P_{3}P_{4}  \cap  P_{19}P_{57} \cap P_{24}P_{28},\\
P_{44}&= P_{15}P_{22} \cap  P_{18}P_{65} \cap  P_{24}P_{28},\\
P_{46}&= P_{3}P_{5} \cap  P_{24}P_{28}\cap P_{51}P_{68},\\
P_{47}&= P_{2}P_{5} \cap  P_{25}P_{38}\cap P_{42}P_{54},\\
P_{49}&=P_{22}P_{23} \cap  P_{19}P_{57} \cap  P_{25}P_{38},\\
P_{50}&= P_{2}P_{4} \cap  P_{18}P_{65}\cap  P_{25}P_{38},\\
P_{55}&= P_{2}P_{5} \cap  P_{36}P_{40}\cap P_{41}P_{52},\\
P_{59}&= P_{3}P_{4} \cap  P_{18}P_{65} \cap  P_{36}P_{40},\\
P_{60}&= P_{25}P_{38} \cap P_{26}P_{66}\cap P_{36}P_{40},\\
P_{61}&= P_{1}P_{6} \cap  P_{36}P_{40}\cap P_{51}P_{68},\\
P_{62}&= P_{1}P_{7} \cap  P_{30}P_{53}\cap P_{42}P_{54},\\
P_{63}&= P_{24}P_{28} \cap  P_{27}P_{56}\cap  P_{30}P_{53},
\end{array}
$$
\end{minipage}
\hskip 20pt
\begin{minipage}{0.45\linewidth}
$$
\begin{array}{ll}
P_{64}&= P_{2}P_{4}\cap  P_{19}P_{57} \cap  P_{30}P_{53},\\
P_{67}&= P_{3}P_{5} \cap  P_{30}P_{53}\cap P_{41}P_{52},\\
P_{69}&= P_{2}P_{4} \cap P_{41}P_{52}\cap  P_{42}P_{54},\\
P_{71}&= P_{1}P_{6}\cap  P_{19}P_{57} \cap  P_{42}P_{54},\\
P_{73}&= P_{3}P_{5}\cap  P_{18}P_{65} \cap  P_{42}P_{54},\\
P_{74}&= P_{2}P_{5}\cap  P_{19}P_{57} \cap  P_{51}P_{68},\\
P_{75}&= P_{1}P_{7}\cap  P_{18}P_{65} \cap  P_{51}P_{68},\\
P_{76}&= P_{22}P_{23}\cap  P_{26}P_{66} \cap  P_{51}P_{68},\\
P_{77}&= P_{3}P_{4}\cap P_{41}P_{52} \cap  P_{51}P_{68},\\
P_{78}&= P_{11}P_{17} \cap  P_{27}P_{56}\cap P_{41}P_{52},\\
P_{79}&= P_{3}P_{5}\cap  P_{19}P_{57} \cap  P_{27}P_{56},\\
P_{81}&= P_{2}P_{5}\cap  P_{18}P_{65} \cap  P_{26}P_{66},\\
P_{82}&= P_{12}P_{17} \cap  P_{26}P_{66}\cap P_{41}P_{52},\\
P_{83}&= P_{3}P_{9} \cap  P_{19}P_{57}\cap P_{41}P_{52},\\
P_{84}&= P_{2}P_{8} \cap  P_{18}P_{65}\cap P_{41}P_{52}.
\end{array}
$$
\end{minipage}
\vskip 5pt

 Removing the factors which lead to the degenerated cases we are left with the equation
\begin{equation*}
\begin{aligned}
C_{24}(a,b):=& a^8b^3 + a^7b^3 + a^6b^4 - 6a^7b^2 + 3a^6b^3 - 6a^6b^2 -   2a^5b^3 + 10a^6b - \\& 6a^5b^2- 2a^4b^3 - a^6 + 12a^5b + 3a^4b^2 - 2a^3b^3 - 6a^5 + 3a^4b + \\
& 6a^3b^2 - a^2b^3 -3a^4 - 13a^3b + 9a^2b^2 - ab^3 + 4a^3 - 12a^2b +  \\& 6ab^2 - b^3 + 5a^2 -3ab +
 2a - b=0.
\end{aligned}
\end{equation*}

This condition is necessary for the construction to terminate successfully in the sense that we obtain exactly $85$ triple points on $24$ lines satisfying the combinatorial properties of the original B\"or\"oczky configuration of $24$ lines.

\subsubsection{Parameter space of configurations $\mathbb{B}_{24}$}
The parametrizing curve
\begin{equation*}
\begin{aligned}
C_{24}(a,b):=& a^8b^3 + a^7b^3 + a^6b^4 - 6a^7b^2 + 3a^6b^3 - 6a^6b^2 -   2a^5b^3 + 10a^6b - \\& 6a^5b^2- 2a^4b^3 - a^6 + 12a^5b + 3a^4b^2 - 2a^3b^3 - 6a^5 + 3a^4b + \\
& 6a^3b^2 - a^2b^3 -3a^4 - 13a^3b + 9a^2b^2 - ab^3 + 4a^3 - 12a^2b +  \\& 6ab^2 - b^3 + 5a^2 -3ab +
 2a - b=0
\end{aligned}
\end{equation*}
 is a curve of genus 5. Therefore, there are no general methods to determine all the rational points of $C_{24}$. The only rational points of height up to $10^5$ are $$(1 : 1: 1),\hspace{0,3cm} (0 : 1 : 0),\hspace{0,3cm} (0: 0 : 1),\hspace{0,3cm} (1 : 0 : 0),\hspace{0,3cm} (-1 : -1 : 1).$$
Simple but tiresome  calculations confirm that all of these points lead to a degene\-ra\-ted configuration, however we were not able to prove that the construction cannot be made over the rational numbers.

\section{Conclusions}

We finish the paper presenting some open problems which we hope to investigate in the future.

\begin{problem}
    Is the configuration $\mathbb{B}_{12}$ the only rational counterexample to the containment $I^{(3)}\subset I^2$ among B\"or\"oczky confi\-gu\-rations?
\end{problem}

 \begin{problem}
    Can we proceed with  the construction of a parameter curve of $\mathbb{B}_{n}$ for any $n$? Can we say something more about the parameter curves?
    Can another choice of points as parameters give a ``nicer'' (easier to investigate) curve $C_{n}$?
 \end{problem}

And finally
 \begin{problem}
    Prove  that all $\mathbb{B}_{n}$ configurations for $n\geq 12$ give counterexamples to the containment $I^{(3)}\subset I^2$.
\end{problem}

\paragraph*{Acknowledgement.}
   These notes originated from discussions during workshops in Lanckorona.\\
   The research of Lampa-Baczy\'nska was partially supported by National Science Centre, Poland, grant 2016/23/N/ST1/01363,  the research of Tutaj-Gasi\'nska was partially supported by National Science Centre, Poland, grant
   2014/15/B/ST1/02197.

   We warmly thank  Maciej Ulas  for helpful conversations and for his constant help with Magma. We also thank Beata Hejmej, Grzegorz Malara, Tomasz Szemberg and Justyna Szpond for valuable suggestions and comments on the text.



\bigskip \small

\bigskip
   \L ucja Farnik,
   Institute of Mathematics, Polish Academy of Sciences, \'Sniadeckich 8,
PL-00-656 Warszawa, Poland
\\
\nopagebreak
   \textit{E-mail address:} \texttt{Lucja.Farnik@gmail.com}

\bigskip
  Jakub Kabat, Magdalena~Lampa-Baczy\'nska,
   Instytut Matematyki UP,
   Podchor\c a\.zych 2,
   PL-30-084 Krak\'ow, Poland
\\
\nopagebreak
  \textit{E-mail address:} \texttt{xxkabat@gmail.com}\\
  \textit{E-mail address:} \texttt{lampa.baczynska@wp.pl}

\bigskip
   Halszka Tutaj-Gasi\'nska,
   Jagiellonian University, Faculty of Mathematics and Computer Science, {\L}ojasiewicza 6, PL-30-348 Krak\'ow, Poland
\\
\nopagebreak
   \textit{E-mail address:} \texttt{Halszka.Tutaj@im.uj.edu.pl}


\end{document}